\documentclass[12pt, twoside, a4paper, fleqn]{article}

\usepackage{latexsym}
\usepackage{amscd}
\usepackage{theorem}
\usepackage{pifont}
\usepackage{mathbbol}
\usepackage{amsfonts}
\usepackage{xspace}
\usepackage{amssymb}
\usepackage{fancyhdr}
\usepackage{mathrsfs}
\usepackage{amsmath}%
\usepackage{graphics}%
\usepackage{graphicx}%
\usepackage{euscript}%
\usepackage{psfrag}%
\usepackage{empheq}%
\newcommand{\cW}{\mathcal{W}}

\newcommand{\bx}{\mathbf{x}}
\newcommand{\by}{\mathbf{y}}
\newcommand{\bb}{\mathbb}
\newcommand{\La}{\Lambda}
\newcommand{\wrt}[1]{\mathrm{d}{#1}}
\newcommand{\vv}{\mathbf{v}}

\newcommand{\R}{{\mathbb R}}

\newcommand{\Z}{{\mathbb Z}}
\newcommand{\N}{{\mathbb N}}
\newcommand{\C}{{\mathbb C}}
\newcommand{\Half}{{\mathbb H}}
\newcommand{\Q}{{\mathbb Q}}

\newcommand{\SL}{{\rm SL}}

\newcommand{\ep}{ \varepsilon }

\newtheorem{theorem}{Theorem}

\newtheorem{proposition}{Proposition}

\theoremstyle{definition}
\newtheorem{definition}[theorem]{Definition}

\theoremstyle{remark}

\newcommand{\dist}{\mathrm{dist}\,}

\newcommand{\be}{\begin{eqnarray*}}
\newcommand{\ee}{\end{eqnarray*}}

\newcommand{\height}{\textrm{\rm H}}
\newcommand{\p}{\bb P}

\newcommand{\pbf}{\mathbf{p}}

\newcommand{\Bad}{\mathbf{Bad}}
\newcommand{\bad}{\mathbf{Bad}}

\newcommand{\est}{\operatorname{EST}}
\newcommand{\EST}{Erd\H{o}s-Sz\"usz-Tur\'an }

\newcommand{\vol}{\operatorname{vol}}

\newcommand{\SO}{\operatorname{SO}}
\newcommand{\Id}{\operatorname{Id}}
\newcommand{\Res}{\operatorname{Res}}
\newcommand{\rar}{\rightarrow}
\newcommand{\cqW}{\mathcal{W}_{p}(\psi, \Q, n)}
\newcommand{\prim}{\operatorname{prim}}

\numberwithin{equation}{section}

\newcommand{\rr}{\mathbf{r}}
\newcommand{\x}{\mathbf{x}}

\begin{document}

\title{Topics in homogeneous dynamics and number theory \bigskip \bigskip }

\author{\textsc{Anish Ghosh} \bigskip \bigskip \\
Tata Institute of Fundamental Research \bigskip \bigskip \\ 
\texttt{email: ghosh@math.tifr.res.in} \bigskip \bigskip }

\date{\today}

\maketitle 

\thispagestyle{empty}

\newpage 

\tableofcontents

\thispagestyle{empty}

\newpage

\section{Introduction} 

This is a survey of some topics at the interface of dynamical systems and number theory, based on lectures delivered at CIRM Luminy, the University of Houston, and IIT Delhi. Specifically, we will be interested in the ergodic theory of group actions on homogeneous spaces and its connections to metric Diophantine approximation. The topics covered in the lectures included the study of the Diophantine approximation of linear forms using dynamics, the study of quadratic forms in particular the famous Oppenheim's conjecture and its variations, as well as lattice point counting using dynamics. At IIT, nondivergence estimates for unipotent flows and Margulis' proof of the Borel Harish-Chandra theorem using the nondivergence estimates were also covered. There are many recent and excellent surveys covering all this material, including but not restricted to \cite{BRV, EW2, E, K0, K2, O} and the books \cite{EW1, WM1}. Rather than reinvent the wheel, I have chosen to present some other recent topics at the interface of Diophantine approximation and homogeneous dynamics in this article. The topics chosen are representative of the lectures but reflect my interests and problems that I have been recently involved with. While some of the lectures were at a more basic level, this article serves as an introduction to more advanced and more recent material. In particular, this is not meant to be a comprehensive survey of this very active and rapidly expanding subject, a shortcoming redressed by the many aforementioned surveys. The hope is that this article will serve as a guide for students with some preparation, e.g. the ones who attended the lectures and point them to further reading and interesting research avenues. Two sections of this article are devoted to results using methods from classical number theory, an indispensable part of the toolbox of anyone interested in Diophantine analysis. 

\subsection{Homogeneous Dynamics}

Let $G$ be a unimodular, locally compact, second countable topological group and $\Gamma$ be a lattice in $G$. The homogeneous space $G/\Gamma$ is equipped with a finite measure which descends from the Haar measure on $G$ and which can therefore be normalised to make $G/\Gamma$ a probability space.  A subgroup $H$ of $G$ acts on the probability space $G/\Gamma$ by translations. The ergodic theory of this action has been extensively studied in recent decades, and is referred to as ``homogeneous dynamics".  For particular choices of $G, H$ and $\Gamma$, the spaces $G/\Gamma$ and $H\backslash G$ parametrise objects of number theoretic interest in many cases and the ergodic theory of the $H$ action on $G/\Gamma$ (resp. the $\Gamma$ action on $H \backslash G$) gives valuable Diophantine information about these objects. Here are some examples:
\begin{enumerate}
\item $G = \SL_{n}(\R)$ and $\Gamma = \SL_{n}(\Z)$. Then $G/\Gamma$ can be identified with the space of unimodular lattices in $\R^n$. The dynamics of diagonal flows on $G/\Gamma$ plays an important role in Diophantine approximation of vectors and linear forms as explained in the next section.

\item Let $n = p+q$ and set $H = \SO(p, q)$. Then the $H$ action on $\SL_{n}(\R)/\SL_{n}(\Z)$ plays an important role in the study of quadratic forms. This is also touched upon in the next section.

\item $G = \SL_{n}(\R) \times \SL_{n}(\Q_p)$ and $\Gamma = \SL_{n}(\Z[1/p])$. Once again, $\Gamma$ is a non-cocompact lattice in $G$ and similar to the example above, one can identify $G/\Gamma$ can be identified with the space of discrete $\Z[1/p]$-modules in $\R^n \times \Q^{n}_p$. Dynamics on this and related spaces plays an important role in $p$-adic Diophantine approximation.\\
 
\item Let $k$ be a degree $d$ number field, $S$ be the set of Archimedean places and $O_k$ be it's ring of integers. Then $O_k$ is a lattice in $k_S :=\R \times \dots \times \R \times \C \times \dots \times \C$ via the Galois embedding. By the Borel Harish-Chandra theorem, $\Gamma = \SL_{2}(O_k)$ is a lattice in $\SL_{n}(k_S) := \prod_{s \in S} \SL_{n}(k_s)$ where $k_s$ denotes the completion of $k$ at the place $s \in S$.  In section \ref{nf}, we will consider flows on the space $G/\Gamma$ and connections to Diophantine approximation in number fields. Such dynamics is intimately connected to geodesic flows on the associated arithmetic orbifold.

\end{enumerate}

\subsection{Diophantine approximation}
Diophantine approximation begins with a theorem due to Dirichlet \cite{Cas} which states that
\begin{theorem}
For any $ x \in \R$ and any $Q > 0$, there exist $ p \in \Z$ and $q \in \N$ such that
$$|q x - p| < 1/Q \text{ and } q \leq Q.$$ 
\end{theorem}

\noindent As a corollary, it follows that for every $x \in \R$, there exist infinitely many $q \in \N$ such that
$$ |q x - p| < 1/q $$
\noindent for some $p \in \Z$. One can also consider Diophantine approximation in higher dimension. Indeed, the role of homogeneous dynamics is brought into sharper focus in higher dimensions as it serves as a replacement for continued fractions, an efficient theory of which is only available in dimension $1$. For instance, the corollary to Dirichlet's theorem in arbitrary dimension reads as follows:
  
\begin{theorem}
For every ${\bf x} \in \mathbb{R}^d$ there exist infinitely many $q \in  \mathbb{Z}$ such that
\begin{equation}\label{dirichlet}
\|q{\bf x} - {\bf p}\| < |q|^{-1/d}.
\end{equation}
for some ${\bf p} \in \mathbb{Z}^d$.
\end{theorem}
Here, $\|~\|$ is the supremum norm. In dimension $d > 1$, there are two possible settings for Diophantine approximation. The above is called the simultaneous setting, one could also consider the \emph{dual} setting where one considers small values of the linear form
$$ |\mathbf{q}\cdot {\bf x} + p| $$ 
for ${\bf q} \in \Z^d$ and $p \in \Z$. The simultaneous and dual settings are related by \emph{transference principles}. We refer the reader to \cite{Cas} for Khintchine's classical transference principles and \cite{CGGMS} for some recent developments involving transference inequalities in the weighted and inhomogeneous settings. Dirichlet's theorem can be proved using the pigeonhole principle (as proved originally by Dirichlet) and also using Minkowski's theorem in the geometry of numbers.

The next major theorem in metric Diophantine approximation seeks to expand the class of approximating functions. Let $\psi$ be a non-increasing function from $\R \to \R_{+} \cup \{0\}$ be given and let $\mathcal{W}_{d}(\psi, \R)$ be the subset of real numbers $x$ for which there exist infinitely many ${\bf q} \in \mathbb{Z}^d$ such that 
\begin{equation}\label{psi1} 
| {\bf q} \cdot {\bf x} + p| < \psi(\|{\bf q}\|^d)
\end{equation}

\noindent for some $p \in \mathbb{Z}$.  Khintchine's theorem (in dual form) characterizes the size of $\mathcal{W}_{d}(\psi, \R)$ in terms of Lebesgue measure.

\begin{theorem}\label{KhinGro}[Khintchine's theorem]
$\mathcal{W}_{d}(\psi, \R)$ has zero or full measure according as 
\begin{equation}\label{grosum}
\sum_{x = 1}^{\infty} \psi(x)
\end{equation}
\noindent converges or diverges. 
\end{theorem}

There exist numbers (and vectors) for which (\ref{dirichlet}) cannot be improved; these are called \emph{badly approximable}. In other words, ${\bf x}$ is called badly approximable if there exists $c:= c({\bf x}) > 0$ such that
\begin{equation}\label{bad}
\|q{\bf x} - {\bf p}\| \geq \frac{c}{|q|^{1/d}}.
\end{equation}

It is well known that badly approximable vectors have zero Lebesgue measure and full Hausdorff dimension (Jarnik \cite{J} for $n=1$ and Schmidt \cite{S1, S3} for arbitrary $n$). In fact, Schmidt showed that they are \emph{winning} for a certain game, a stronger and more versatile property than having full Hausdorff dimension. We refer the reader to Dani's article \cite{D2} in this volume for an introduction to Schmidt's game.\\ 

On the opposite end of the spectrum to badly approximable vectors, are singular vectors. A vector $\x \in \R^d$ is said to be \emph{singular} if for every $ \ep > 0 $ there exists $ N_0$ with the following property: for each $ N \ge N_0$, there exist $\pbf \in \Z^d$, $ q \in \N $ so that
\begin{equation}\label{singclassical}
\| q{\bf x} - {\bf p} \| < \frac{\ep}{  N^{\frac1d}} \qquad {\rm and } \qquad q < N \, .
\end{equation}
In other words, $ \x $ is singular if Dirichlet's Theorem can be improved by an arbitrarily small constant factor $\ep>0$. In the case $d=1$, Khintchine \cite{kh} showed that a real number is singular if and only if it is rational. Moreover, it was shown by Davenport $\&$ Schmidt \cite{DavSch}  that the set of singular vectors has zero Lebesgue measure.\\ 

We now move from linear forms to quadratic forms and briefly discuss Oppenheim's conjecture. Let $n \geq 3$ and let $Q$ be a non degenerate indefinite quadratic form in $n$ variables and assume that $Q$ is not proportional to a form with rational coefficients. It was a conjecture of Oppenheim from the 1920's and a celebrated theorem of Margulis \cite{Margulis} that under these conditions $Q(\Z^n)$ is dense in $\R$. Oppenheim's conjecture is false for binary quadratic forms, a counterexample can be constructed using badly approximable numbers. More precisely, the quadratic form $Q(x, y) = y^2 - \theta^2 x^2$ where $\theta$ is a quadratic irrational with $\theta^2$ irrational provides a counterexample. For this and more, we refer the reader to Borel's survey \cite{Borel}.\\

How does the Diophantine approximation of vectors and linear and quadratic forms relate to dynamics of subgroup actions on $G/\Gamma$ or lattice actions on $H \backslash G$?  Let $G = \SL_{n+1}(\R)$ and $\Gamma = \SL_{n+1}(\Z)$, then $G/\Gamma$ and can be naturally identified with the space $\Omega_{n+1}$ of unimodular, i.e. covolume $1$ lattices in $\R^{n+1}$. Namely, $G$ acts transitively on $\Omega_{n+1}$ by multiplication and the stabilizer of the lattice $\Z^{n+1}$ is $\SL_{n+1}(\Z)$. The space $G/\Gamma$ is a non-compact, finite volume space and Mahler's compactness criterion describes the compact subsets of $G/\Gamma$. Diophantine approximation of vectors in $\R^n$ can be modelled using dynamics of subgroup actions on $\Omega_{n+1}$.   Given a vector ${\bf x} \in \R^n$ we consider the unimodular lattice 
$$ \Lambda_{\bf x} := \begin{pmatrix} 1 & {\bf x}\\ 0 & \Id \end{pmatrix} \Z^{n+1} \in \Omega_{n+1}.$$
Further, let
$$g_t :=  \begin{pmatrix} e^t & 0 \\ 0 & e^{-t} \end{pmatrix} \in G.$$
Then we have the following two propositions connecting Diophantine properties of ${\bf x}$ with the dynamics of the $g_t$ action on $G/\Gamma$ due to Dani \cite{D1}. The first concerns badly approximable vectors.
\begin{proposition}
A vector ${\bf x} \in \R^n$ is badly approximable if and only if $\{g_t\Lambda_{\bf x}:~t~ > 0 \}$ is bounded in $G/\Gamma$.
\end{proposition}

And the second concerns singular vectors.

\begin{proposition}
A vector ${\bf x} \in \R^n$ is singular if and only if $\{g_t\Lambda_{\bf x}:~t~ > 0 \}$ is divergent in $G/\Gamma$.
\end{proposition}

Kleinbock and Margulis \cite{KM2} have proved a more general version of the ``Dani correspondence"and have provided a dynamical proof of Khintchine's theorem using exponential mixing of the $g_t$ action on $G/\Gamma$. This is closely related to the \emph{shrinking target problem} for group actions on homogeneous spaces.  In \cite{Su}, Sullivan established the following folklore theorem.  Let $V = \mathbb{H}^{d+1}/\Gamma$ be a hyperbolic manifold where $\Gamma$ is a discrete subgroup of hyperbolic isometries which is not co-compact, and let $\dist v(t)$ denote the distance from a fixed point in $V$ of the point achieved after traveling a time $t$ along the geodesic with initial direction $v$.

\begin{theorem}\label{K-S-theorem}\cite{Su}
For all $x \in V$, and almost every $v \in T_x V$,
$$ \limsup_{t \to \infty} \frac{\dist v(t)}{\log t}  = 1/d.$$

\end{theorem}

Kleinbock and Margulis generalized Sullivan's logarithm law to locally symmetric spaces. In fact, both the logarithm law and Khintchine's theorem are manifestations of a $0-1$ Borel-Cantelli type law for diagonal flows on homogeneous spaces. 
This scheme was subsequently carried out in the positive characteristic setting in \cite{AGP} (see also \cite{GR}), and in the setting of geodesic orbits on the frame bundle of finite volume non-compact hyperbolic manifolds in \cite{ABG}. We refer the reader to \cite{GK1} for an introduction to the shrinking target problem, including a more general formulation of the logarithm law for (not necessarily one-parameter") subgroups, as well as a list of references.\\

 How does Oppenheim's conjecture relate to dynamics?  Let $G = \SL_{n}(\R)$, $H = \SO(p, q)$ and $\Gamma = \SL_{n}(\Z)$. The space $H\backslash G$ can be identified with the space of quadratic forms of signature $(p, q)$ in $n = p + q$ variables. The relevant dynamics here is that of the $\Gamma$ action on $H\backslash G$ (or, dually, the $H$ action on $G/\Gamma$). More precisely, the following Proposition implies Oppenheim's conjecture. 
\begin{proposition}
Any $H$ orbit on $G/\Gamma$ is either closed and carries an $H^{\circ}$ invariant measure or is dense.
\end{proposition}  
The above result was proved by G. Margulis \cite{Margulis} for ternary quadratic forms, thus settling Oppenheim's conjecture. The main point is that, under the conditions above, $\SO(p, q)$ is generated by unipotent one-parameter subgroups. The Proposition above is an instance of the conjectures of Raghunathan and Dani on orbit closures and invariant measures for actions of such groups on homogeneous spaces. These conjectures were settled by M. Ratner. We refer the reader to \cite{WM1} for details on this beautiful subject.
Recently, the $\Gamma$ action on $H\backslash G$ has been used to study \emph{effective} versions on Oppenheim's conjecture.  We will not elaborate on this theme in this survey, referring the reader instead to \cite{GGN5} for details and to \cite{GGN3, GGN4} for the more general study of effective density of lattice orbits on homogeneous varieties and \cite{GGN1, GGN2} for the related problem of intrinsic Diophantine approximation on varieties. See also \cite{GK2} for a related question on quadratic forms studied originally by Bourgain \cite{B}.\\

\subsection{Structure}
The rest of the article is divided into three sections. The next section focuses on Dirichlet's theorem and considers two different aspects - probabilistic and geometric, of the problem of \emph{distribution of approximates}. The tools in this section involve equidistribution of flows and limiting distributions for flows on the space of lattices. Section \ref{nf} considers Diophantine approximation in number fields. The classical theorems like Dirichlet's theorem and Khintchine's theorem can be generalised to number fields. We pay particular attention to vectors which are badly approximable by rationals from a number field and the associated dynamics of diagonal flows on arithmetic orbifolds. The last two sections discuss Diophantine approximation in two diverse settings: that of projective space and hyperbolic space. In section \ref{projective}, we present a projective version of Khintchine's theorem and the more general Duffin-Schaeffer conjecture and in section \ref{hyperbolic}, we discuss Diophantine approximation by orbits of Fuchsian and Kleinian groups on the boundary of hyperbolic space. The techniques used in the last two sections are classical in nature.

\section*{Acknowledgements} This survey grow out of lectures delivered at CIRM Luminy, at the University of Houston and at IIT Delhi. I am grateful to the organisers of each of the three events for inviting me and for their hospitality. Special thanks to Jayadev Athreya, Alan Haynes and Riddhi Shah. I would also like to thank the editors of this volume, Anima Nagar, Riddhi Shah and Shrihari Sridharan. This work was supported by a grant from the Indo-French Centre for the Promotion of Advanced Research; a Department of Science and Technology, Government of India Swarnajayanti fellowship; a MATRICS grant from the Science and Engineering Research Board; and the Benoziyo Endowment Fund for the Advancement of Science at the
Weizmann Institute. I gratefully acknowledge the hospitality of the Technion and the Weizmann Institute.

\section{On the distribution of approximates}
In this section we describe some recent results on the distribution of approximates in Dirichlet's theorem. First, we describe a probabilistic distribution problem, originally due to Erd\H{o}s, Sz\"usz and Tur\'an and developed in a homogeneous context, in \cite{AG}. In the next subsection, we consider the geometry of the approximates and describe a spiraling equidistribution of approximates proved in \cite{AGT1}.

\subsection{The EST distribution}

\noindent In 1958, Erd\H{o}s, Sz\"usz and Tur\'an  ~\cite{EST} introduced a problem in probabilistic Diophantine approximation: what is the probability $f(N, A, c)$ that a point $\alpha$ chosen from the uniform distribution on $[0,1]$ has a solution $\frac p q \in \Q$ to the inequality 
\begin{equation}\label{DM}
 \left| \alpha - \frac p q \right| \le \frac {A}{q^2} 
 \end{equation} 
 with the constraint that the denominator $q$ lies in $[N, cN]$? Here $A >0, c>1$ are fixed positive parameters, and $N$ is a parameter which goes to infinity. The above inequality is of course a close variant of the inequality in Dirichlet's theorem. In particular, we know that $A =1$ admits infinitely many solutions and that (by Hurwitz's theorem), $A = \frac{1}{\sqrt{5}}$ is the best allowable constant which admits infinitely many solutions for all $\alpha$.  Given $A, c, N$, let $\est(A,c,N)(\alpha)$ be the number of solutions $p/q \in \Q$ with $\gcd(p, q) = 1$ to (\ref{DM}). Letting $\alpha \in [0, 1]$ be a uniform random variable yields an integer-valued random variable $\est(A,c,N)$, with  $$P(\est(A,c,N) = k ) = m\left(\alpha \in [0, 1]: \mbox{ there are exactly } k \mbox{ solutions to (\ref{DM})} \right),$$ where $m$ is Lebesgue measure on $[0,1]$. Then, the \EST question is the existence of the limit $$\lim_{N \rightarrow \infty} P(\est(A,c, N)>0).$$

\noindent  Kesten~\cite{Kesten} considered a modifed version of this problem, he defined the sequence of random variables $K(A, N)$ as the number of solutions to
\begin{equation}\label{DM1} \left| \alpha q- p \right| \le \frac {A}{N}, 1 \le q \le N,\end{equation} where $\alpha$ is a uniform $[0,1]$ random variable. That is, $$P(K(A,N) = k ) = m\left(\alpha \in [0, 1]: \mbox{ there are exactly } k \mbox{ solutions to (\ref{DM1})} \right).$$

The Kesten distribution was studied by Marklof [Theorem 4.4 in \cite{Marklof}], thought at the time he was not aware of Kesten's question. In \cite{AG}, it was shown that the the limiting distributions of the random variables $\est(A,c, N)$ and $K(A, N)$ exist as $N \rightarrow \infty$. In fact, they can be viewed as the probability of a random unimodular lattice intersecting a certain fixed region. Let $\mu_2$ denote the Haar probability measure on $\Omega_2$, the space of unimodular lattices in $\R^2$. Given $\La \in X_2, \La = g \Z^2$, let $\La_{\prim}$ be the set of primitive vectors in $\La$. It should be emphasized that we were unaware of Marklof's work, and the Kesten part of the Theorem below merely reproves a part of Marklof's theorem referred to above. 

\begin{theorem}\label{theorem:est:lattice}(\cite{AG} Theorem $1.1$) The limiting distribution of the random variables $\est(A,c,N)$ and $K(A, N)$ exist and denoting the random variables with these limiting distributions as $\est(A,c)$ and $K(A)$, we have \begin{equation}\label{eq:est:dist}P(\est(A,c) = k) = \mu_2(\La \in X_2: \#(\La_{\prim} \cap H_{A,c}) = k ),\end{equation} and 
\begin{equation}\label{eq:kesten:dist}P(K(A) = k) = \mu_2(\La \in X_2: \#(\La_{\prim} \cap R_{A}) = k )
\end{equation} where
\begin{equation}\label{def:region}
H_{A, c} = \{(x,y) \in \R^2~:~xy \leq A, 1 \leq y \leq c\},
\end{equation}
and
\begin{equation}\label{def:region:rectangle}
R_{A} = \{(x,y) \in \R^2~:~|x| \leq A, 0 \leq y \leq 1\}.
\end{equation}

\end{theorem}

In fact, the setting in \cite{AG} is abstract and axiomatic, allowing for a great deal of flexibility. The philosophy of \emph{equivariant point processes}, introduced in this context in \cite{AG}, allows us to obtain the existence of these limiting distributions in higher dimensions, for linear forms, for points on smooth curves as well as in the setting  of the set of holonomy vectors of saddle connections on translation surfaces. For details as well as a more detailed history of the problem, the reader is referred to \cite{AG}. Further applications of equivariant point processes are explored in a forthcoming monograph of Athreya and Ghosh.

\subsection{Spiraling of approximates}

In this section, we describe some results from \cite{AGT1} (see also \cite{AGT2}) where the geometric study of the distribution of approximates in Dirichlet's theorem was initiated. We consider a vector ${\bf x} \in \R^d$ and form, as before, form the associated unimodular lattice in $\R^{d+1}$
\begin{equation}\nonumber
\Lambda_{\bf x} := \begin{pmatrix} \Id_{d} & {\bf x}\\0 & 1 \end{pmatrix} \bb Z^{d+1} = \left\{\begin{pmatrix} q{\bf x} - {\bf p}\\ q \end{pmatrix} ~:~{\bf p} \in \bb Z^{d}, q \in \bb Z\right\}.
\end{equation}

\noindent Then we can view the approximates $({\bf p}, q)$ of ${\bf x}$ appearing in (\ref{dirichlet}) as points of the lattice $\Lambda_{\bf x}$ in the region
\begin{equation}\label{cone}
R := \left\{{\bf v} = \begin{pmatrix}{\bf v}_1\\ v_2 \end{pmatrix} \in \bb R^{d} \times \bb R~:~\|{\bf v}_1\||v_2|^{1/d} \leq 1\right\}.  
\end{equation} 

The goal is to understand the geometry of the set of approximates $\La_{\bf x} \cap R$. To do so, consider the following sets:

\begin{equation}\label{defR1}
R_{\epsilon, T} := \left\{ {\bf v} \in  R~:~ \epsilon T \le v_2 \le T \right\}  
\end{equation}

\noindent and, for a measurable subset $A$ of $\bb S^{d-1}$ with zero measure boundary,

\begin{equation}\label{defR2}
R_{A, \epsilon, T} := \left\{ {\bf v} \in R_{\epsilon, T}~:~ \frac{{\bf v}_1}{\|{\bf v}_1\|} \in A \right\}.
\end{equation} 

\noindent For a unimodular lattice $\La$, define 

$$N(\Lambda, \epsilon, T) = \#\{\Lambda \cap R_{\epsilon, T}\}$$ and 

$$N(\Lambda, A, \epsilon, T) = \#\{\Lambda \cap R_{A, \epsilon, T}\}$$

\noindent Let $dk$ denote Haar measure on $K = \SO_{d+1}(\bb R)$, and let $X_{d+1} = \SL_{d+1}(\bb R)/\SL_{d+1}(\bb Z)$ denote the space of unimodular lattices in $\R^{d+1}$. The following equidistribution theorem is proved in \cite{AGT1}.

\begin{theorem}\label{main}
For every $\Lambda \in X_{d+1}$, $A \subset \bb S^{d-1}$ as above, and for every $\epsilon > 0$, 
\begin{equation}\label{main-1}
 \lim_{T \rightarrow \infty} \frac{\int_{K} N(k^{-1}\Lambda, A, \epsilon, T)~\wrt k}{\int_{K} N(k^{-1}\Lambda, \epsilon, T)~\wrt k}= \vol(A). 
\end{equation}
\end{theorem}
\medskip
\noindent In other words, for any lattice $\La$, on average over the set of directions $\vv$, the set of approximates satisfying Dirichlet's theorem in the direction $\vv$ equidistributes in the set of directions $\bb S^{d-1}$ in the orthogonal complement to $\vv$. The proof of Theorem \ref{main} depends on an equidistribution result for Siegel transforms which is likely to have other applications. 

One can also fix the vertical, and instead average the counting functions over a range of heights, $T$, to obtain a result for almost every lattice $\La$, with respect to the probability measure $\mu$ on $X_{d+1}$ induced by Haar measure on $\SL_{d+1} (\R)$. 

\begin{theorem}\label{mainbirkhoff} Fix $A \subset \bb S^{d-1}$ as above. For $\mu$-almost every $\Lambda \in X_{d+1}$ and for every $\epsilon > 0$,
\begin{equation}\label{main-birkhoff}
 \lim_{S \rightarrow \infty} \frac{\int_{0}^{S} N(\Lambda, A, \epsilon, e^t)~\wrt t}{\int_{0}^S N(\Lambda, \epsilon, e^t)~\wrt t}= \vol(A). \end{equation}
\end{theorem}
\medskip
\noindent That is, if we average the number of approximates in the region $A$ over a range of heights and similarly average the total number of approximates, we have an almost everywhere equidistribution. 

On the other hand, there are examples of lattices $\La$ and directions $\vv$ for which (non-averaged) equidistribution does not hold. The following is proved in \cite{AGT1}.

\begin{theorem}\label{theomaineg}
Let $d \geq 1$.  There exists a lattice $\Lambda \in \SL_{d+1}(\bb R) / \SL_{d+1}(\bb Z)$, a set $A \subset \bb S^{d-1}$ with zero measure boundary, and a sequence $\{T_n\}$ for which \begin{eqnarray*} \label{eqnmaineg}
\lim_{n \rightarrow \infty}  \frac {N(\Lambda, A, \epsilon, T_n)}{N(\Lambda, \epsilon, T_n)} \neq \vol(A)
\end{eqnarray*}
for every $1 > \epsilon\geq0$.
\end{theorem}

\noindent For $d =1$, note that $\bb S^0:=\{-1, 1\}$ and we define $\vol(\{-1\})=\vol(\{1\}) = 1/2$. \\\\

In \cite{KSW}, the study of weighted spiraling was taken up and several interesting spiraling and equidistribution results were obtained. Subsequently, in \cite{AlamG}, a study of spiraling and equidistribution in number fields was undertaken. Finally, we note that the paper \cite{BMPS} has interesting results on the distribution of approximates for badly approximable numbers.

\section{Diophantine approximation in number fields}\label{nf}
In this section, we describe some recent advances in Diophantine approximation in number fields. Let $k$ be a number field of degree $d$ over $\Q$, $O_k$ its ring of integers, and $S$ be the set of field embeddings $\sigma: k\hookrightarrow\R$. Then we have $|S|=d.$ We will be interested in Diophantine approximation by rationals from a \emph{fixed} number field. 
Analogues of Dirichlet's theorem in this setting have been established by several authors (cf. \cite{S3}, \cite{Bur}, \cite{Q}) using appropriate adaptations of the geometry of numbers. Here is Proposition $2.1$ from \cite{KL}, proved using a result from \cite{Bur}, here we are approximating ${\bf x} \in \R^d$ by rationals in $k$.

\begin{proposition}
There exists a constant and for every $Q>0$, there exists $p \in O_k, q \in O_k \backslash \{0\}$ with 
$$\|\sigma(q) \cdot {\bf x} + \sigma(p)\| \leq CQ^{-1}\text{  and } \|\sigma(q)\| \leq Q.$$
\end{proposition}

The corollary of Dirichlet's theorem also holds for number fields. Here is Theorem $2.3$ in \cite{KL}.

\begin{theorem}
There is a constant $C = C_k > 0$ depending only on $k$, such that for every $x \in k_S$, there are infinitely many $p, q\in O_K$ with $q\neq 0$ satisfying:
$$\|\sigma(q) \cdot {\bf x} + \sigma(p)\| \leq C\|\sigma(q)\|^{-1}.$$
\end{theorem}

We now want to define and discuss the properties of badly approximable vectors. We define the more general \emph{weighted} badly approximable vectors and describe results in \cite{AGGL}, we therefore follow the notation from that paper. Let $\rr\in \R^d$ be a real vector with $r_{\sigma}\ge 0$ for $\sigma\in S$ and $\sum_{\sigma\in S}r_{\sigma}=1$. Set
$$S_1=\{\sigma\in S:r_{\sigma}>0\}, \text{ and } S_2=S\setminus S_1.$$
Assume $|S_1|=d_1, |S_2|=d_2$.  Choose and fix $\omega\in S$ with $r_{\omega}=r$, where
$$r=\max_{\sigma\in S}r_{\sigma}. $$
 Define a weighted norm, called the \emph{$\rr$-norm}, on $\prod_{\sigma\in S}\R$ by
$$\|\bx\|_{\rr}=\max_{\sigma\in S_1}|x_{\sigma}|^{\frac{1}{r_{\sigma}}}.$$
\begin{definition}
Say a vector $\bx=(x_{\sigma})_{\sigma\in S} \in \prod_{\sigma\in S}\R$ is \emph{$(k,\rr)$-badly approximable } if
\begin{equation*}
\inf_{\substack{q\in O_k\setminus \{0\}\\p\in O_k }} \max \left\{\max_{\sigma\in S_1} \|q\|_{\rr}^{r_{\sigma}}|\sigma(q)x_{\sigma}+\sigma(p)|,
\max_{\sigma\in S_2} \max\{|\sigma(q)x_{\sigma}+\sigma(p)|,|\sigma(q)|\} \right\}>0.
\end{equation*}
The set of $(k,\rr)$-badly approximable vectors is denoted as $\Bad(k,\rr)$.
\end{definition}

The definition of $(k,\rr)$-badly approximable vector is the weighted case of $k$-badly approximable vector introduced in \cite{EGL}. Weighted badly approximable (by rational) vectors in $\R^n$ are the subject of Schmidt's conjecture, now a theorem of Badziahin, Pollington and Velani \cite{BPV}. In \cite{AGGL}, a number field analogue of Schmidt's conjecture is proved, this was previously known in some special cases \cite{EGL}. The existence of $k$-badly approximable vectors was established in \cite{B} and \cite{H}.
We note in passing that  a real number is badly approximable  if and only if its
partial fraction coefficients are bounded. In \cite{Hines1}, this characterization is established for complex numbers and in \cite{Hines2}, examples of badly approximable vectors in the number field setting have been constructed.\\

A variant of Schmidt's game, called the hyperplane potential game was introduced in \cite{FSU} and defines a class of subsets of $\R^{d}$
called \emph{hyperplane potential winning} (\emph{HPW} for short) sets.

The hyperplane potential game involves two parameters $\beta\in(0,1)$ and $\gamma>0$.  Bob starts the game
 by choosing a closed ball $B_0\subset \R^{d}$ of radius $\rho_0$. In the $i$-th turn,
Bob chooses a closed ball $B_i$ of radius $\rho_i$, and then Alice chooses a countable family of hyperplane
neighborhoods $\{L_{i,k}^{(\delta_{i,k})}: k\in \N\}$ such that
\begin{equation*}
\sum_{k=1}^\infty \delta_{i,k}^\gamma\le(\beta \rho_{i})^\gamma.
\end{equation*}
Then in the $(i+1)$-th turn, Bob chooses a closed ball $B_{i+1}\subset B_i$ of radius $\rho_{i+1}\ge\beta \rho_{i}$. By
this process there is a nested sequence of closed balls
$$B_0\supseteq B_1\supseteq B_2\supseteq \cdots.$$
We say a subset $S\subset \R^{d}$ is \emph{$(\beta,\gamma)$-hyperplane potential winning}
(\emph{$(\beta,\gamma)$-HPW} for short) if no matter how Bob plays, Alice can ensure that
$$\bigcap_{i=0}^\infty B_i\cap\Big(S\cup\bigcup_{i=0}^\infty\bigcup_{k=1}^\infty
L_{i,k}^{(\delta_{i,k})}\Big)\ne\emptyset.$$ We say $S$ is \emph{hyperplane potential winning}
(HPW) if it is $(\beta,\gamma)$-HPW for any $\beta\in(0,1)$ and $\gamma>0$.

Set
\begin{equation*}
\theta: k\to \prod_{\sigma\in S} \R, \quad \theta(p)=(\sigma(p))_{\sigma\in S}.
\end{equation*}
Let $\Res_{k/\Q}$ denote Weil's restriction of scalar's functor. It is well known that the group $\Res_{k/\Q}\SL_2(\Z)$
is a lattice in  $\Res_{k/\Q}\SL_2(\R)$. The latter coincides with the product of $d$ copies of $\SL_2(\R)$. We set
$$G=\Res_{k/\Q}\SL_2(\R)=\prod_{\sigma\in S}\SL_2(\R), \quad \Gamma=\Res_{k/\Q}\SL_2(\Z).$$
It follows from the definition that the subgroup $\Res_{k/\Q}\SL_2(\Z)$ coincides with the subgroup $\theta(\SL_2(O_k))$, where $\theta$ is the map defined by $\theta(g)=(\sigma(g))_{\sigma\in S}$.
The following is a special case of the main Theorem in \cite{AGGL}.

\begin{proposition}\label{main prop}
Let $\rr\in \R^d$ be a real vector with $r_{\sigma}\ge 0$ for $\sigma\in S$ and $\sum_{\sigma\in S}r_{\sigma}=1$, set
\begin{equation}
g_{\rr}(t):=  \left(\begin{pmatrix}e^{r_{\sigma}t}&0\\0&e^{-r_{\sigma}t}\end{pmatrix}\right)_{\sigma\in S}
\end{equation}
and $F_{\rr}^+=\{g_{\rr}(t):t\ge 0\}$, then the set
\begin{equation*}
E(F_{\rr}^+):=\{x\in G/\Gamma: F_{\rr}^+x \text{ is bounded }\}
\end{equation*}
is HPW.
\end{proposition}

As before, $\Bad(k,\rr)$ corresponds to bounded orbits for certain flows on homogeneous spaces. Namely,  we have the following correspondence (Proposition $3.4$ in \cite{AGGL}) between $(k,\rr)$-badly approximable vector and bounded $F_{\rr}^+$ trajectories, i.e. a number field version of Dani's correspondence. The proof is more involved than the $\Q$ case. 
\begin{proposition}\label{dani}
 A vector $\bx=(x_{\sigma})_{\sigma\in S}$ is $(k,\rr)$-badly approximable if and only if the trajectory $F_{\rr}^+u(\bx)\Gamma$ is bounded in $G/\Gamma$.
In other words,
\begin{equation}\label{e:dani}
\Bad(K,\rr)= u^{-1}\big(\pi^{-1}(E(F_{\rr}^+))\cap H\big),
\end{equation}
where $\pi$ denotes the projection $G\rightarrow G/\Gamma$.
\end{proposition}

In view of the number field Dani correspondence above, this implies a number field version of Schmidt's conjecture, in other words.

\begin{theorem}
$\Bad(k,\rr)$ is HPW.
\end{theorem}

In fact, the sets above are winning for the hyperplane absolute game introduced in \cite{BFKRW}. This is because it is proved in \cite[Theorem C.8]{FSU} that A subset $S$ of $\R^{d}$ is HPW if and only if it is HAW.

\section{A Projective Duffin Schaeffer Theorem}\label{projective}
In this section, we describe a recent \emph{projective} variation of metric Diophantine approximation originally studied by Choi and Vaaler \cite{CV} and developed further by the author and Haynes \cite{GH}. Projective metric Diophantine approximation aims to quantify the density of $\p^{n-1}(k)$ in $\p^{n-1}(k_v)$ where $k$ is a number field and $k_v$ is a completion of $k$. For non-zero vectors $\bx, \by \in k^{n}_v$ we define following \cite{CV}
\begin{equation}\label{def-metric}
\delta_{v}(\bx, \by) := \frac{|\bx \wedge \by|_v}{|\bx|_v |\by|_v}.
\end{equation}
Then $\delta_v$ defines a metric on $\p^{n-1}(k_v)$ which induces the usual quotient topology (\cite{Rumely}). We define the height of a point $\bx\in \p^{n-1}(k)$ by
\begin{equation}\label{height-def}
\height(\bx) := \prod_{v}|\bx|_v,
\end{equation}
and we note that this is well defined over projective space because of the product formula. The following is a projective version of Dirichlet's theorem due to Choi and Vaaler \cite{CV}.
\begin{theorem}\label{Choi-Vaaler}
Let $\bx \in \p^{n-1}(k_v)$, let $\tau \in k_v$ with $|\tau|_v \geq 1$. Then there exists $\by \in \p^{n-1}(k)$ such that
\begin{enumerate}
\item $\height(\by) \leq c_{k}(n)|\tau|_{v}^{n-1}$,\text{ and}\\
\item $\delta_{v}(\bx, \by) \leq c_{k}(n)(|\tau|_v\height(\by))^{-1}$.
\end{enumerate}
\end{theorem}
\noindent Here
\begin{align}
c_{k}(n) = 2|\Delta_{k}|^{1/2d} \prod_{v \arrowvert \infty} r_{v}(n)^{d_v/d},\label{constant}
\end{align}
$\Delta_k$ is the discriminant of $k$, and
\begin{align*}
r_{v}(n) = \left\{
\begin{array}{rl}  \pi^{-1/2}\Gamma(\frac{n}{2} + 1)^{1/n} & \text{if } v \text{ is real},\\\\ (2\pi)^{-1/2}\Gamma(n + 1)^{1/2n} & \text{if } v \text{ is complex}.
\end{array} \right.
\end{align*}

In \cite{GH}, a projective analogue of Khintchine's theorem and more generally, the Duffin Schaeffer conjecture, were proved. In order to state the results in loc. cit. we first briefly recall some probability measures on $\p^{n-1}(k_v)$, originally defined and studied by Choi \cite{Choi}. If $v$ is an infinite place then $\beta_v^n$ is the usual $n$-fold Lebesgue measure on $\R^n$ or $2^n$ times Lebesgue measure on $\C^n$, while if $v$ is a finite place then $\beta_v^n$ is the $n-$fold Haar measure normalized so that
$$ \beta_v(O_v) = \|\mathcal{D}_v\|_{v}^{d_v/2}, $$
where $O_v$ is the ring of integers of $k_v$ and $\mathcal{D}_v$ is the local different of $k$ at $v$. Let $\phi :k_v^n\setminus \{{\bf 0}\}\rar\p^{n-1}(k_v)$ be the quotient map and define the $\sigma$-algebra $\mathcal{M}$ of measurable sets in $\p^{n-1}(k_v)$ to be the collection of sets $M\subseteq\p^{n-1}(k_v)$ such that $\phi^{-1}(M)$ lies in the $\sigma-$algebra of Borel sets in $k_v^n$. This is in fact the $\sigma$-algebra of Borel sets in $\p^{n-1}(k_v)$. One then defines measures $\mu_v$ on $(\p^{n-1}(k_v),\mathcal{M})$ by
\begin{equation}\label{def:measure}
\mu_v(M)=\frac{\beta_v^n\left(\phi^{-1}(M)\cap B({\bf 0},1)\right)}{\beta_v^n\left(B({\bf 0},1)\right)}.
\end{equation}

Given $\psi : \R_{+} \cup \{0\} \to \R_{+} \cup \{0\}$ let $\cW$ be the set of $\bx \in \p^{n-1}(\Q_v)$ for which there exist infinitely many $\by \in \p^{n-1}(\Q)$ such that
$$ \delta_{v}(\bx, \by) \leq \psi(\height(\by)).$$



\noindent Then it is a straightforward consequence of the Borel-Cantelli lemma that $\mu_{p}(\cW) = 0$ whenever
\begin{equation}\label{sum}
\sum_{q=1}^{\infty} q^{n - 1}\psi(q)^{(n - 1)}
\end{equation}
\noindent converges. In \cite{GH}, the projective $p$-adic version of the Duffin-Schaeffer conjecture is established in all dimensions greater than $1$.

\begin{theorem}\label{thm2}
Assume that $p$ is a finite place, that $n > 2$, and that $\psi (q)=0$ whenever $p|q$. Then $\mu_p(\cqW) = 1$ whenever (\ref{sum}) diverges.
\end{theorem}

\noindent In fact, more can be proved if monotonicity is assumed. The second result in \cite{GH} is the complete (i.e. allowing arbitrary primes and dimensions) projective version of Khintchine's theorem.

\begin{theorem}\label{thm1}
Assume that $\psi$ is decreasing and let $p$ be a (finite or infinite) place of $\Q$. Then $\mu_{p}(\cqW)=1$ whenever (\ref{sum})
diverges.
\end{theorem}

\noindent Recently, in \cite{HH}, the authors have studied badly approximable vectors in the projective setting. In particular, they showed that badly approximable vectors have full Hausdorff dimension.

\section{The hyperbolic picture}\label{hyperbolic}
Consider the action of $\SL_{2}(\R)$ on the hyperbolic upper half plane $\Half^2$ by M\"{o}bius transformations. This action extends to the boundary and there are close connections between Diophantine approximation and the study of dense orbits of discrete subgroups of $\SL_{2}(\R)$ on the boundary. For example, the orbit of $\infty$ under the action of $\SL_{2}(\Z)$ is precisely the set of rational numbers; one might therefore seek a more general quantitative theory of  the approximation of limit points of a fixed Kleinian group by points in the orbit (under the group) of a distinguished limit point $y$. In this section, we briefly review important work by Patterson \cite{Pat2} and then discuss some recent work in this direction carried out in \cite{BGSV}.

Let $G$ denote\footnote{This notation, mainstream in the Kleinian groups literature, is at odds with the notation in previous sections where $G$ was the ambient Lie group and $\Gamma$ a lattice in $G$.} a nonelementary, geometrically finite Kleinian group acting on the unit ball model $(B^{d+1}, \rho)$ of $(d+1)$--dimensional hyperbolic space with metric $ \rho$ derived from the differential $ d \rho = 2 | d \x | /(1-|\x |^2 ) $. Thus, $G$ is a discrete subgroup of $ \mbox{M\"ob}(B^{d+1})$, the group of orientation-preserving M\"obius transformations of the unit ball $B^{d+1}$. Since $G$ is nonelementary, the limit set $\Lambda$ of $G$ is uncountable. The group $G$ is said to be of the first kind (such a group is a lattice) if $\Lambda = S^d $ and of the second kind otherwise. Let $ \delta $ denote the Hausdorff dimension of $\Lambda$.  It is well known that $\delta $ is equal to the exponent of convergence of the group. For  $g \in G$ set $L_g := |g^\prime(0)|^{-1}$, where $|g^\prime(0)| = 1-|g(0)|^2 $ is the (Euclidean) conformal dilation of $g$ at the origin. It can be checked that $L_g \leq e^{\rho(0,g(0))} \leq 4 L_g$.  The following two Dirichlet-type theorems were first established by Patterson \cite[Section 7: Theorems 1 \& 2]{Pat2} for finitely generated Fuchsian groups, but can be generalized to higher dimensions.
\begin{theorem}\label{pat1}
Let $G$ be a nonelementary, geometrically finite Kleinian group containing parabolic elements and let $P$ be a complete set of inequivalent parabolic fixed points of $G$. Then there is a constant $c> 0 $ with the following property: for each $\xi \in \Lambda$, $ N > 1$, there exist $p \in P$, $g \in G$ so that
\[
| \xi - g(p) | \le \frac{c}{\sqrt{L_g N}}
\qquad and \qquad L_g \le N \, .
\]
\end{theorem}

\begin{theorem}
Let $G$ be a nonelementary, geometrically finite Kleinian group without parabolic elements and let $\{\eta,\eta'\}$ be the pair of fixed points of a hyperbolic element of $G$. Then there is a constant $c> 0 $ with the following property: for all $\xi \in \Lambda$, $ N > 1$, there exist $y \in \{\eta,\eta'\}$, $g \in G$ so that
\[
| \xi - g(y) | \le \frac{c}{ N} \qquad and \qquad L_g \le N \, .
\]
\end{theorem}

As mentioned earlier, in the context of $\SL_{2}(\Z)$ and $\Half^2$, Theorem \ref{pat1} reduces to Dirichlet's Theorem.   

In \cite{BGSV}, the notion of singular limit points within the hyperbolic space setup was introduced. Let $G$ be a Kleinian group and let $Y$ be a complete set $P$ of inequivalent parabolic fixed points of $G$ if the group has parabolic elements; otherwise let $Y$ be the pair $\{\eta,\eta'\}$ of fixed points of a hyperbolic element of $G$. A point $\xi \in \Lambda $ is said to be \emph{singular} if for every $ \ep > 0 $ there exists $ N_0$ with the following property: for each $ N \ge N_0$, there exist $y \in Y$, $g \in G$ so that
\begin{equation}\label{singhyper}
| \xi - g(y) | \ < \ \left\{
\begin{array}{ll}
\frac{\ep}{\sqrt{L_g N}} \quad & {\rm if}
\;\;\; Y = P \; \\ [4ex]
\frac{\ep}{ N } \quad & {\rm if} \;\;\; Y = \{\eta,\eta'\} \;
\end{array}
\right.
\qquad {\rm and } \qquad
L_g < N \ . \end{equation}

\noindent In \cite{BGSV}, it was shown that the hyperbolic ``singular'' theory is, irrespective of the dimension of the hyperbolic space, similar to the one-dimensional classical theory.

\begin{theorem} \label{singthm}
Let $G$ be a nonelementary, geometrically finite Kleinian group, and let $Y$ be as above. Then a point $ \xi \in \Lambda $ is singular if and only if $\xi \in G(Y) := \{ g(y) : g \in G, y \in Y \}$.
\end{theorem}

In \cite{Pat2}, convergence and divergence Khintchine type theorems were proved. In \cite{BGSV}, versions of Khintchine's theorem for \emph{proper} subsets of the limit set were investigated.  Let $K$ be a subset of the limit set $\Lambda$ which supports a nonatomic probability measure $\mu$.  Given $\alpha > 0$, the measure $\mu$ supported on $K$ is said to be {\em weakly absolutely $\alpha$-decaying}  if there exist strictly positive constants $ C, r_0 $ such that for all
$\ep > 0$ we have
\[
\mu\big(B(x,\ep r) \big) \ \leq \ C \, \ep^{\alpha} \mu\big(B(x,r)\big) \hspace{7mm} \forall \ x \in K
\hspace{5mm} \forall \ r < r_0 \ .
\]
For sets supporting such measures, the following result was proved in \cite{BGSV}.

\begin{theorem} \label{mainext}
Let $G$ be a nonelementary, geometrically finite Kleinian group and let $y$ be a parabolic fixed point of $G$, if there are any, and a hyperbolic fixed point otherwise. Fix $\alpha > 0$, and let $K$ be a compact subset of $\Lambda$ equipped with a weakly absolutely $\alpha$-decaying measure $\mu$. Then
\begin{equation} \label{ineqmain2}
\mu( K \cap W_{y}(\psi) ) = 0\hspace{6mm} {\rm if \ } \hspace{6mm} \sum_{r=1}^{\infty}
r^{\alpha -1 } \psi(r)^{\alpha} \ < \ \infty \ .
\end{equation}
\end{theorem}

We now discuss the analogue of badly approximable vectors. The set
\[
\bad_y := \left\{ \xi \in \Lambda : \ \exists \ \ c(\xi) > 0 {\rm \ such \ that \ } | \xi - g(y) | > c(\xi)/ L_g \ \ \forall \ \
g \in G \right\} \, ,
\]
can be considered to be the hyperbolic analogue of badly approximable numbers and is of measure zero. Nevertheless, it is a large set. The following theorem was first established by Patterson \cite[Section 10]{Pat2} for finitely generated Fuchsian groups of the first kind. As before, $y$ is taken to be a parabolic fixed point of $G$ if the group has parabolic elements and a hyperbolic fixed point of $G$ otherwise.

\vspace{3ex}

\begin{theorem}
Let $G$ be a nonelementary, geometrically finite Kleinian group and let $y$ be a parabolic fixed point of $G$, if there are any, and a hyperbolic fixed point otherwise. Then
\[
\dim \bad_y = \dim \Lambda \, .
\]
\end{theorem}

Let $K$ be a subset of the limit set $\Lambda$ which supports a nonatomic probability measure $\mu$ as before. We assume that the measure $\mu$ supported on $K$ is {\em Ahlfors $\delta$-regular} for some $\delta > 0$; that is, that there exist constants $ C > 0$ and $r_0 $ such that
\[
C^{-1} \, r^{\delta}\ \le \ \mu\big(B(x, r) \big) \ \le \ C \, r^{\delta} \hspace{7mm} \forall \ x \in K
\hspace{5mm} \forall \ r < r_0 \ .
\]
\noindent  Sets supporting such measures are referred to as Ahlfors $\delta$-regular and it is a well known fact that
 $$
 \dim K =  \delta  \, .
 $$
 For Ahlfors $\delta$-regular  subsets of the limit set  the following result was proved in \cite{BGSV}.

\begin{theorem} \label{mainbad}
Let $G$ be a nonelementary, geometrically finite Kleinian group and let $y$ be a parabolic fixed point of $G$, if there are any, and a hyperbolic fixed point otherwise. Let $K$ be a compact, Ahlfors $\delta$-regular subset of $\Lambda$. Then
\begin{equation} \label{ineqmain3}
\dim \left( K \cap \bad_y \right) = \dim K \ . \end{equation}
\end{theorem}

These results were motivated by results on Diophantine approximation on manifolds, and indeed, constitute an hyperbolic analogue of the theory. We refer to \cite{BGSV} for details.

\end{document}